# Changement de base pour les foncteurs $Tor$

21 novembre 2018

Mathieu Zimmermann
Institut de Recherche Mathématique Avancée de Strasbourg
Université Louis Pasteur et CNRS
7, rue René Descartes
67 084 Strasbourg CEDEX France
zimmerma@math.u-strasbg.fr

Dans [5], Pirashvili démontre un théorème de changement de base pour les foncteurs dérivés $Ext$ qui s'exprime sous la forme $Ext^*_\mathcal{B}(La, b) \cong Ext^*_\mathcal{C}(a, Rb)$ où $L$ et $R$ sont des foncteurs adjoints, exacts et $L$ (resp. $R$) préserve les projectifs (resp. les injectifs). Le but de ce papier est de donner la notion analogue à celle de foncteur adjoint pour avoir un théorème de changement de base pour les foncteurs dérivés $Tor$. De plus, on définira une classe de petites catégories sur lesquelles on peut appliquer le résultat et qui permettra une réinterprétation du résultat de T. Pirashvili et B. Richter donné dans [6].

## 0 Rappels et notations

Soit $K$ un anneau commutatif. On renvoie à [4] pour les notions de la théorie des catégories. Pour toute petite catégorie $\mathcal{B}$, on note $\mathcal{B}$-$Mod$ la catégorie des foncteurs de $\mathcal{B}$ dans la catégorie des $K$-modules. Les morphismes entre deux $\mathcal{B}$-modules $M$ et $N$ sont les transformations naturels et on se permet l'abus de notation $Hom_\mathcal{B}(M, N) := Hom_{\mathcal{B}-Mod}(M, N)$. Enfin, si $M$ un $\mathcal{B}^{op}$-module et $N$ un $\mathcal{B}$-module, on définit le produit tensoriel au dessus de $\mathcal{B}$ par :

$$M \otimes_{\mathcal{B}} N := \left( \bigoplus_{X \in Ob(\mathcal{B})} M(X) \otimes_K N(X) \right) \Big/ \sim$$

où $\sim$ est la relation d'équivalence donnée par :

$$M(\Phi)m \otimes_K n \sim m \otimes_K N(\Phi)n$$

où $\Phi \in Hom_{\mathcal{B}}(X,Y)$, $m \in M(Y)$, $n \in N(X)$. Notons que $M \otimes_{\mathcal{B}} N = N \otimes_{\mathcal{B}^{op}} M$. On construit alors les foncteurs dérivés $Tor_n^{\mathcal{B}}(M,N) := H_n(p_* \otimes_{\mathcal{B}} N) = H_n(M \otimes_{\mathcal{B}} q_*)$ où $p_*$ est une résolution projective de $M$ et $q_*$ est une résolution projective de $N$. Dans la suite, $\mathcal{C}$ et $\mathcal{D}$ désigneront deux petites catégories. Si $\mathcal{C}$ est une sous-catégorie de $\mathcal{D}$, alors $\mathcal{O}$ sera, suivant le cas, le foncteur oubli des $\mathcal{D}$-modules dans les $\mathcal{C}$-modules ou le foncteur oubli des $\mathcal{D}^{op}$-modules dans les $\mathcal{C}^{op}$-modules.

# 1    Théorème de changement de base

**1.1 Définition : pseudo-adjonction.** Soient

$$\begin{array}{rccl} L: & \mathcal{C}^{op}\text{-}Mod & \longrightarrow & \mathcal{D}^{op}\text{-}Mod, \\ R: & \mathcal{D}\text{-}Mod & \longrightarrow & \mathcal{C}\text{-}Mod, \end{array}$$

deux foncteurs. Ils donnent naissance à deux bifoncteurs, covariants en chaque variable :

$$L- \otimes_{\mathcal{D}} -\ \text{ et }\ - \otimes_{\mathcal{C}} R- : \mathcal{C}^{op}\text{-}Mod \times \mathcal{D}\text{-}Mod \longrightarrow K\text{-}Mod.$$

S'il existe un isomorphisme naturel entre ces deux bifoncteurs, on dit qu'ils sont *pseudo-adjoints*. Le foncteur $L$ sera dit *pseudo-adjoint à gauche* de $R$ et le foncteur $R$ sera dit *pseudo-adjoint à droite* de $L$ par rapport aux catégories $\mathcal{C}$ et $\mathcal{D}$ (par rapport aux catégories $\mathcal{C}^{op}$ et $\mathcal{D}^{op}$, le foncteur $L$ serait le foncteur pseudo-adjoint à droite).

**1.2 Théorème de changement de base pour les foncteurs $Tor$.** *Soient $L : \mathcal{C}^{op}\text{-}Mod \longrightarrow \mathcal{D}^{op}\text{-}Mod$ et $R : \mathcal{D}\text{-}Mod \longrightarrow \mathcal{C}\text{-}Mod$ deux foncteurs pseudo-adjoints. Si l'une des conditions suivantes est remplie :*

1. *$R$ est exact et envoie projectif sur projectif,*
2. *$L$ est exact et envoie projectif sur projectif,*

alors, pour tout $\mathcal{C}^{op}$-module $a$ et pour tout $\mathcal{D}$-module $b$, on a l'isomorphisme naturel suivant :
$$Tor_*^{\mathcal{C}}(a, Rb) \cong Tor_*^{\mathcal{D}}(La, b).$$

*Démonstration :* Supposons que la condition (1) soit remplie et soit $q_*$ une résolution projective de $b$, alors :
$$\begin{aligned} Tor_*^{\mathcal{D}}(La, b) &= H_*(La \otimes_{\mathcal{D}} q_*) \\ &\cong H_*(a \otimes_{\mathcal{C}} Rq_*) \end{aligned}$$

car les foncteurs $L$ et $R$ sont pseudo-adjoints. Comme $R$ est exact et qu'il envoie projectif sur projectif, $Rp_*$ est une résolution projective de $Rb$ et donc $H_*(a \otimes_{\mathcal{C}} Rp_*) = Tor_*^{\mathcal{C}}(a, Rb)$. Le deuxième cas se fait de manière identique. □

## 2 Adjonction et pseudo-adjonction

Le but de cette section est de donner la relation entre l'adjonction et la pseudo-adjonction qui, dans le cas où on se limite aux $K$-modules libres de rang fini, se trouvent être deux notions très proches.
Rappelons que tout $K$-module $V$ définit un foncteur contravariant :
$$\begin{aligned} Hom_K(-, V) : \mathcal{C}\text{-}Mod &\longrightarrow \mathcal{C}^{op}\text{-}Mod \\ M &\longmapsto Hom_K(M, V). \end{aligned}$$

**2.1 Définition.** Dans le cas où $V = K$, nous noterons :
$$\begin{aligned} * : \mathcal{C}\text{-}Mod &\longrightarrow \mathcal{C}^{op}\text{-}Mod \\ M &\longmapsto M^* := Hom_K(M, K) : X \mapsto M(X)^*. \end{aligned}$$

Si $F : \mathcal{C}\text{-}Mod \to \mathcal{D}\text{-}Mod$ est un foncteur, nous pouvons définir un nouveau foncteur $\overline{F} := * \circ F \circ *$ :
$$\begin{aligned} \overline{F} : \mathcal{C}^{op}\text{-}Mod &\longrightarrow \mathcal{D}^{op}\text{-}Mod \\ M &\longmapsto (F(M^*))^*. \end{aligned}$$

**2.2 Lemme.** *Soit $M$ un $\mathcal{C}$-module. Le foncteur $- \otimes_{\mathcal{C}} M : \mathcal{C}^{op}\text{-}Mod \to K\text{-}Mod$ est adjoint à gauche au foncteur $Hom_K(M, -) : K\text{-}Mod \to \mathcal{C}^{op}\text{-}Mod$, i.e. si $N$ est un $\mathcal{C}^{op}$-module et $V$ est un $K$-module on a :*
$$Hom_K(N \otimes_{\mathcal{C}} M, V) \cong Hom_{\mathcal{C}^{op}}(N, Hom_K(M, V)).$$

*En particulier, en prenant $V = K$ :*
$$(N \otimes_{\mathcal{C}} M)^* \cong Hom_{\mathcal{C}^{op}}(N, M^*).$$

*Démonstration :* Soit $\Phi : N \to Hom_K(M, V)$ un morphisme de $\mathcal{C}^{op}$-modules. On peut alors définir une application $K$-linéaire :

$$\begin{aligned}
\tilde{\Phi} : N \otimes M & \longrightarrow V \\
n \otimes m & \longmapsto (\Phi(n))(m)
\end{aligned}$$

Cette valeur reste inchangée si on prend un autre élément dans la même classe d'équivalence de $N \otimes_{\mathcal{C}} M$. En effet, si $g \in \mathcal{C}$ :

$$\begin{aligned}
\tilde{\Phi}(ng \otimes m) & = (\Phi(ng))(m) \\
& = (\Phi(n) \circ g)(m) \quad \text{car } \Phi \in Hom_{\mathcal{C}^{op}}(N, Hom_K(M, V)) \\
& = \Phi(n)(gm) \\
& = \tilde{\Phi}(n \otimes gm)
\end{aligned}$$

et $\tilde{\Phi} \in Hom_K(N \otimes_{\mathcal{C}} M, V)$.

Soit maintenant $\Psi \in Hom_K(N \otimes_{\mathcal{C}} M, V)$. On peut alors définir une application linéaire :

$$\begin{aligned}
\Psi' : N & \longrightarrow Hom_K(M, V) \\
n & \longmapsto \Psi(n \otimes_{\mathcal{C}} -) : \begin{array}{rcl} M & \to & V \\ m & \mapsto & \Psi(n \otimes_{\mathcal{C}} m). \end{array}
\end{aligned}$$

De la même manière que dans le premier cas, on vérifie qu'il y a compatibilité avec la structure de $\mathcal{C}^{op}$-module.

Finalement, comme ces applications sont l'inverse l'une de l'autre, on a bien prouvé l'isomorphisme. $\square$

**2.3 Théorème.** *Notons $\mathcal{C}\text{-}Mod'$ la sous-catégorie pleine de $\mathcal{C}\text{-}Mod'$ constituée des foncteurs à valeurs dans les $K$-modules libres de rang fini. Soient $L : \mathcal{C}^{op}\text{-}Mod' \to \mathcal{D}^{op}\text{-}Mod'$ et $R : \mathcal{D}\text{-}Mod' \to \mathcal{C}\text{-}Mod'$ deux foncteurs pseudo-adjoints. Dans ce cas, le foncteur $L$ est adjoint à gauche du foncteur $\overline{R} = * \circ R \circ *$.*

*Démonstration :* Soient $N$ un $\mathcal{C}^{op}$-module libre de rang fini et $M$ un $\mathcal{D}^{op}$-module libre de rang fini. Alors, comme $M^{**} = M$ :

$$\begin{aligned}
Hom_{\mathcal{D}^{op}}(L(N), M) & \cong (L(N) \otimes_{\mathcal{D}} M^*)^* && \text{d'après le lemme 2.2,} \\
& \cong (N \otimes_{\mathcal{C}} R(M^*))^* && \text{grâce à la pseudo-adjonction,} \\
& \cong Hom_{\mathcal{C}^{op}}(N, R(M^*)^*) && \text{d'après le lemme 2.2,} \\
& = Hom_{\mathcal{C}^{op}}(N, \overline{R}(M)).
\end{aligned}$$

$\square$

**Remarque :** si on reprend la démonstration, il suffit de demander que $M$

soit libre de rang fini pour qu'on ait un isomorphisme. sectionCatégories croisées

On va donner ici des exemples pour lesquels on peut construire des paires de foncteurs pseudo-adjoints.

**2.4 Définition.** Soient $\mathcal{B}$ une petite catégorie et $\mathcal{C}$ et $\mathcal{D}$ deux sous-catégories ayant les mêmes objets que $\mathcal{B}$. Si pour tout $\Phi \in Hom_\mathcal{B}(X,Z)$ il existe un unique $Y \in Ob(\mathcal{B})$, un unique $\Psi \in Hom_\mathcal{C}(Y,Z)$ et un unique $f \in Hom_\mathcal{D}(X,Y)$ tels que $\Phi = \Psi \circ f$, on dira que $\mathcal{B}$ est une *catégorie croisée* par $\mathcal{C}$ et $\mathcal{D}$ et on la notera $\mathcal{B} = \mathcal{C} \bowtie \mathcal{D}$.

Des exemples de ce type de catégories sont donnés par :
- la catégorie simpliciale qui est croisée par la sous-catégorie engendrée par ses injections et la sous-catégorie engendrée par ses surjections (voir par exemple [3]),
- les catégories simpliciales croisées (voir [2]),
- les produits semi-directs de groupes,
- le groupe symétrique $\Sigma_n = \mathbb{Z}/n\mathbb{Z} \bowtie \Sigma_{n-1}$.

**2.5 Actions des catégories.** Soient $\mathcal{B} = \mathcal{C} \bowtie \mathcal{D}$ et $f \in Hom_\mathcal{D}(Y,Z)$ et $\Psi \in Hom_\mathcal{C}(X,Y)$. Comme $f \circ \Psi \in Hom_\mathcal{B}(X,Z)$, $f$ et $\Psi$ déterminent uniquement $Y' \in Ob(\mathcal{B})$, $\Psi^*(f) \in Hom_\mathcal{D}(X,Y')$ et $f_*(\Psi) \in Hom_\mathcal{C}(Y',Z)$ tels que

$$f \circ \Psi = f_*(\psi) \circ \Psi^*(f),$$

ce qui se lit sur le diagramme suivant :

$$\begin{array}{ccc} & \Psi & \\ X & \longrightarrow & Y \\ \Psi^*(f) \downarrow & & \downarrow f \\ Y' & \longrightarrow & Z \\ & f_*(\Psi) & \end{array}$$

Il existe alors des applications :

$$f_* : \coprod_{X \in Ob(\mathcal{B})} Hom_\mathcal{C}(X,Y) \longrightarrow \coprod_{Y' \in Ob(\mathcal{B})} Hom_\mathcal{C}(Y',Z),$$

$$\Psi^* : \coprod_{Z \in Ob(\mathcal{B})} Hom_\mathcal{C}(Y,Z) \longrightarrow \coprod_{Y' \in Ob(\mathcal{B})} Hom_\mathcal{C}(X,Y').$$

**2.6 Lemme.** *Soient $\Psi_1$, $\Psi_2$, $\Psi$ des morphismes de $\mathcal{C}$ et $f^1$, $f^2$ et $f$ des morphismes de $\mathcal{D}$. Si $\mathcal{B} = \mathcal{C} \bowtie \mathcal{D}$, on aura, à chaque fois que cela fait sens pour la composition, les relations suivantes :*

- $(f^1 \circ f^2)_*(\Psi) = f^1_*(f^2_*(\Psi))$
- $(\Psi_1 \circ \Psi_2)^*(f) = \Psi_2^*(\Psi_1^*(f))$
- $f_*(\Psi_1 \circ \Psi_2) = f_*(\Psi_1) \circ (\Psi_1^*(f))_*(\Psi_2)$
- $\Psi^*(f^1 \circ f^2) = (f^2_*(\Psi))^*(f^1) \circ \Psi^*(f^2)$
- $id_X^*(f) = f$ et $f_*(id_X) = id_Z$
- $id_{Y*}(\Psi) = \Psi$ et $\Psi^*(id_Y) = id_X$

*Démonstration :* Démontrons par exemple la première et la quatrième assertions.

$$\begin{aligned} f^1 \circ f^2 \circ \Psi &= f^1 \circ f^2_*(\Psi) \circ \Psi^*(f^2), \\ &= f^1_*(f^2_*(\Psi)) \circ (f^2_*(\Psi))^*(f^1) \circ \Psi^*(f^2). \end{aligned}$$

Comme $f^1 \circ f^2 \circ \Psi = (f^1 \circ f^2)_*(\Psi) \circ \Psi^*(f^1 \circ f^2)$, on procède par identification. Les deuxième et troisième assertions se démontrent de façon similaire, les deux dernières sont immédiates. □

**2.7 Construction du foncteur pseudo-libre.** On va montrer que le foncteur oubli $\mathcal{O} : \mathcal{B}\text{-}Mod \to \mathcal{C}\text{-}Mod$ admet un pseudo-adjoint à gauche. On pose :

$$\begin{aligned} \mathcal{L}_\mathcal{D} : \mathcal{C}^{op}\text{-}Mod &\longrightarrow \mathcal{B}^{op}\text{-}Mod \\ M &\longmapsto \bigoplus_{A \,\in\, Ob(\mathcal{B})} M(A) \otimes K[Hom_\mathcal{D}(-, A)] \end{aligned}$$

ou, de manière équivalente :

$$\begin{aligned} \mathcal{L}_\mathcal{D}(M) : \mathcal{B}^{op} &\longrightarrow K\text{-}Mod \\ X &\longmapsto \bigoplus_{A \,\in\, Ob(\mathcal{B})} M(A) \otimes K[Hom_\mathcal{D}(X, A)]. \end{aligned}$$

Si $\Phi = \Psi \circ f$ est la décomposition canonique, on définit son action sur $\mathcal{L}_\mathcal{D}(M)$ par :

$$(m \otimes g)\Phi = mg_*(\Psi) \otimes (\Psi^*(g) \circ f).$$

En utilisant le lemme précédent, on démontre qu'on a bien une structure de $\mathcal{B}^{op}$-module. On peut remarquer que ce foncteur est en plus un foncteur libre, *i.e.* le foncteur adjoint à gauche au foncteur d'oubli de $\mathcal{B}^{op}\text{-}Mod$ dans $\mathcal{C}^{op}\text{-}Mod$.

**2.8 Proposition.** *Soit $\mathcal{B} = \mathcal{C} \bowtie \mathcal{D}$ une catégorie croisée. Alors le foncteur oubli $\mathcal{O} : \mathcal{B}\text{-}Mod \longrightarrow \mathcal{C}\text{-}Mod$ admet un pseudo-adjoint à gauche $\mathcal{L}_\mathcal{D}$. De plus, $\mathcal{L}_\mathcal{D}$ est exact et préserve les projectifs.*

*Démonstration :* Regardons d'abord quelles sont les classes d'équivalences dans $\mathcal{L}_\mathcal{D} M \otimes_\mathcal{B} N$. Soit $M$ (resp. $N$) un $\mathcal{C}^{op}$-module (resp. un $\mathcal{B}$-module). Soient $m \in M_A$, $n \in N_X$ et $g \in Hom_\mathcal{D}(X, A)$. Comme

$$\begin{aligned}(m \otimes g) \otimes n &= (m \otimes id_A)g \otimes n \\ &\sim (m \otimes id_A) \otimes gn \text{ dans } \mathcal{L}_\mathcal{D} M \otimes_\mathcal{B} N,\end{aligned}$$

chaque classe d'équivalence est représentée par un élément de la forme $(m \otimes id_X) \otimes n$. On n'a plus à considérer l'action des automorphismes et il suffit de calculer les classes d'équivalence du type suivant :

$$(m \otimes id)\Psi \otimes n \quad \sim \quad (m \otimes id) \otimes \Psi n$$
$$\iff$$
$$m \otimes n\Psi \quad \sim \quad m \otimes \Psi n$$

où $\Psi$ est un morphisme de $\mathcal{C}$ et on a bien une bijection entre $\mathcal{L}_\mathcal{D}(M) \otimes_\mathcal{B} N$ et $M \otimes_\mathcal{C} \mathcal{O}(N)$.

Pour prouver la naturalité de cette bijection, il suffit de remarquer que les diagrammes suivants sont commutatifs :

$$\begin{array}{ccc} \mathcal{L}_\mathcal{D}(M_1) \otimes_\mathcal{B} N & \longrightarrow & M_1 \otimes_\mathcal{C} \mathcal{O}(N) \\ {\scriptstyle \mathcal{L}_\mathcal{D}(f) \otimes_\mathcal{B} N} \downarrow & & \downarrow {\scriptstyle f \otimes_\mathcal{C} \mathcal{O}(N)} \\ \mathcal{L}_\mathcal{D}(M_2) \otimes_\mathcal{B} N & \longrightarrow & M_2 \otimes_\mathcal{C} \mathcal{O}(N) \end{array}$$

$$\begin{array}{ccc} \mathcal{L}_\mathcal{D}(M) \otimes_\mathcal{B} N_1 & \longrightarrow & M \otimes_\mathcal{C} \mathcal{O}(N_1) \\ {\scriptstyle \mathcal{L}_\mathcal{D}(M) \otimes_\mathcal{B} g} \downarrow & & \downarrow {\scriptstyle M \otimes_\mathcal{C} \mathcal{O}(g)} \\ \mathcal{L}_\mathcal{D}(M) \otimes_\mathcal{B} N_2 & \longrightarrow & M \otimes_\mathcal{C} \mathcal{O}(N_2). \end{array}$$

Il est immédiat que $\mathcal{L}_\mathcal{D}$ est exact car, dans la construction, on ne fait que tensoriser avec un module libre. Soit maintenant $p$ un $\mathcal{C}^{op}$-module projectif. Pour voir que $\mathcal{L}_\mathcal{D}$ préserve les projectifs, il suffit de voir qu'il y a une injection canonique :

$$\begin{aligned} p &\longrightarrow \mathcal{O}(\mathcal{L}_\mathcal{D}(p)) \\ x &\longmapsto x \otimes id \end{aligned}$$

Supposons qu'on ait le diagramme suivant dans $\mathcal{B}^{op}\text{-}Mod$ :

$$\begin{array}{ccc} & & \mathcal{L}_\mathcal{D}(p) \\ & & \downarrow f \\ b & \longrightarrow & a \\ & \pi & \end{array}$$

On sait qu'il existe un $h : p \subset \mathcal{O}(\mathcal{L}_\mathcal{D}(p)) \longrightarrow \mathcal{O}(b)$ tel que $f_{|p} = \pi \circ h$. Il ne reste alors qu'une unique extension $h' : \mathcal{L}_\mathcal{D}(p) \to a$ qui fasse commuter le diagramme, ce qui prouve que $\mathcal{L}_\mathcal{D}$ préserve les projectifs (et que $\mathcal{L}_\mathcal{D}$ est le foncteur libre de la catégorie des $\mathcal{C}^{op}$-modules dans la catégorie des $\mathcal{B}^{op}$-module). □

**2.9 Corollaire.** *Soit $\mathcal{B} = \mathcal{C} \bowtie \mathcal{D}$ une catégorie croisée. Pour tout $\mathcal{C}^{op}$-module $a$ et pour tout $\mathcal{B}$-module $b$, on a l'isomorphisme naturel suivant :*

$$Tor_*^\mathcal{B}(\mathcal{L}_\mathcal{D}(a), b) \cong Tor_*^\mathcal{C}(a, \mathcal{O}(b))$$

*Démonstration :* C'est une application immédiate du théorème 1.2 et de la proposition précédente. □

# 3   Homologie simpliciale et homologie cyclique

Dans [6] sont construits deux foncteurs $b$ et $\overline{b}$ tels que pour tout $\mathcal{F}(as)$-module $M$ (resp. pour tout $\Gamma(as)$-module $N$), on ait des isomorphismes :

$$\begin{aligned} Tor_*^{\Gamma(as)}(\overline{b}, M) &\cong Tor_*^{\Delta^{op}}(K, \mathcal{O}(M)) \\ Tor_*^{\mathcal{F}(as)}(b, N) &\cong Tor_*^{\Delta C^{op}}(K, \mathcal{O}(N)), \end{aligned}$$

où $K$ est le $\Delta^{op}$-module (resp. $\Delta C^{op}$-module) trivial qui envoie chaque objet sur $K$ et chaque morphisme sur l'identité. On verra dans cette partie que ces isomorphismes proviennent de la construction plus générale décrite précédemment et qu'on pourra réinterpréter les foncteurs $b$ et $\overline{b}$ en tant qu'image des modules triviaux par les foncteurs pseudo-adjoints à gauche des foncteurs oublis.

**3.1 Rappels et notations.** (voir [2], [3]). A partir de maintenant, les catégories auront toutes pour objets les $[n] = \{0, \ldots, n\}$ avec $n \in \mathbb{N}$.
La catégorie simpliciale $\Delta$ est la catégorie ayant pour morphismes les applications ensemblistes croissantes au sens large. Elle est engendrée par les morphismes $\sigma_i^n : [n+1] \to [n]$, $0 \leqslant i \leqslant n$, qui envoient $i$ et $i+1$ sur $i$ et les morphismes $\delta_i^n : [n-1] \to [n]$, $0 \leqslant i \leqslant n$, dont l'image ne contient pas $i$. S'il n'y a pas de confusion possible, on n'écrira pas l'exposant. La catégorie opposée $\Delta^{op}$ sera engendrée par les morphismes $d_i^n$ (resp. $s_i^n$) correspondants aux morphismes $\delta_i^n$ (resp. $\sigma_i^n$).
La catégorie $\mathcal{F}(as)$ est la catégorie dont les morphismes sont les applications ensemblistes avec un ordre total sur les images réciproques. La catégorie $\Gamma(as)$ (noté $\Delta^{op}S'$ dans [3]) est la sous-catégorie de $\mathcal{F}(as)$ dont les morphismes envoient 0 sur 0.

Notons $\Sigma_\bullet$ (resp. $\Sigma_{\bullet+1}$) le groupoïde qui vérifie $Hom([n],[m]) = \emptyset$ et $Hom([n],[n])$ est le groupe symétrique $\Sigma_n$ (resp. $\Sigma_{n+1}$). Nous allons considérer deux catégories simpliciales croisées introduites dans [2]. La catégorie $\Delta S = \Delta \bowtie \Sigma_{\bullet+1}$, dans laquelle tout morphisme $\Phi \in Hom_{\Delta S}([n],[m])$ s'écrit de manière unique sous la forme $\Phi = \Psi \circ \sigma$ où $\Psi \in Hom_\Delta([n],[m])$ et $\sigma \in \Sigma_{n+1}$. De même, la sous-catégorie $\Delta C$ de $\Delta S$ est la catégorie dans laquelle chaque morphisme $\Phi \in Hom_{\Delta C}([n],[m])$ s'écrit de manière unique $\Phi = \Psi \circ \sigma$ où $\Psi \in Hom_\Delta([n],[m])$ et $\sigma = t_n^k \in \mathbb{Z}/(n+1)\mathbb{Z}$ où $t_n$ est un générateur du groupe $\mathbb{Z}/(n+1)\mathbb{Z}$. Il est montré dans [6] que $\mathcal{F}(as)$ est isomorphe à $\Delta S$. Rappelons (voir [3] chapitre 6) qu'il existe une inclusion $\Delta C^{op} \hookrightarrow \Delta S$ qui est donnée par :

$$\begin{aligned}
d_i &\mapsto \begin{cases} \sigma_i & \text{si } 0 \leqslant i < n \\ \sigma_0 \circ (n, 0, \ldots, n-1) & \text{si } i = n \end{cases} \\
s_i &\mapsto \delta_i \\
t &\mapsto (n, 0, \ldots, n-1)
\end{aligned}$$

qui donne un isomorphisme entre $\Delta C$ et $\Delta C^{op}$. Enfin, le $\Delta C^{op}$-module trivial $K$ sera le foncteur qui associera $K$ à tout objet et l'application linéaire identité à chaque morphisme. De la même manière, $K$ désignera aussi le $\Delta^{op}$-module trivial.

**3.2 Proposition.** *La catégorie $\mathcal{F}(as) = \Delta S$ est une catégorie croisée par les catégories $\Delta C^{op}$ et $\Sigma_\bullet$, i.e. $\mathcal{F}(as) = \Delta C^{op} \bowtie \Sigma_\bullet$.*

*La catégorie $\Gamma(as)$ est une catégorie croisée par les catégories $\Delta^{op}$ et $\Sigma_\bullet$, i.e. $\Gamma(as) = \Delta^{op} \bowtie \Sigma_\bullet$.*

*Démonstration :* On sait que $\Delta C^{op} \cong \Delta C \subset \Delta S$. Soit $\Phi \in Hom_{\Delta S}([n],[m])$. Comme $\Delta S = \Delta \bowtie \Sigma_{\bullet+1}$, $\Phi$ peut s'écrire de manière unique sous la forme $\Phi = \Psi \circ \sigma$ où $\Psi \in Hom_\Delta([n],[m])$ et $\sigma \in \Sigma_{n+1}$. Or $\sigma$ peut à son tour s'écrire de manière unique $\sigma = t^k \circ \sigma'$ où $t$ est l'image de l'opérateur cyclique de $\Delta C$ et $\sigma' \in \Sigma_n$ (en voyant $\Sigma_n \subset \Sigma_{n+1}$ comme étant le sous-groupe laissant 0 invariant). Finalement, $\Phi \in Hom_{\Delta S}([n],[m])$ peut s'écrire de manière unique :

$$\begin{aligned}
\Phi &= \Psi & \circ & \; \sigma \\
&= \underbrace{\Psi \circ t^k} & \circ & \; \sigma' \\
&= \Psi' & \circ & \; \sigma'
\end{aligned}$$

où $\Psi' \in Hom_{\Delta C}([n],[m]) \cong Hom_{\Delta C^{op}}([n],[m])$ et $\sigma' \in \Sigma_n$, ce qui prouve la première assertion.

Pour vérifier la seconde assertion, remarquons que :

$$\Phi = \Psi' \circ \sigma' \in \Gamma(as) \Leftrightarrow \Psi' \in \Delta^{op}.$$

$\square$

**3.3 Les foncteurs $b$ et $\overline{b}$.** Le $\mathcal{F}(as)^{op}$-module $b$ est défini dans [6] comme étant le conoyau de l'application :

$$\begin{aligned} d : K[Hom_{\mathcal{F}(as)}(-,[1])] &\longrightarrow K[Hom_{\mathcal{F}(as)}(-,[0])] = K[\Sigma_{-+1}] \\ \Phi &\longmapsto (d_0 - d_1) \circ \Phi \end{aligned}$$

où $d_0 = \delta_0$ et $d_1 = \delta_0 \circ t$. Soit $\Phi \in Hom_{\mathcal{F}(as)}([n],[1])$ et $\Phi^{-1}(0) = (i_1 < \ldots < i_m)$ et $\Phi^{-1}(1) = (j_1 < \ldots < j_p)$ avec $m+p = n+1$. Alors, $(d_0 \circ \Phi)^{-1}(0) = (i_1 < \ldots < i_m < j_1 < \ldots < j_p)$ et $(d_1 \circ \Phi)^{-1}(0) = (j_1 < \ldots < j_p < i_1 < \ldots < j_m)$. Les classes d'équivalence sont donc en bijection avec les ordres cycliques sur $[n]$ (*i.e.* un ordre total à permutation cyclique près).

De même, $\overline{b}$ est défini comme étant le conoyau de l'application

$$\overline{d} : K[Hom_{\Gamma(as)}(-,[1])] \longrightarrow K[Hom_{\Gamma(as)}(-,[0])] = K[\Sigma_{-+1}]$$

qui est la restriction de l'application précédente. La description de $\overline{b}$ est identique.

**3.4 Proposition.** *Le $\mathcal{F}(as)^{op}$-module $b$ (resp. le $\Gamma(as)^{op}$-module $\overline{b}$) est l'image par le foncteur $\mathcal{L}_{\Sigma_\bullet}$ du $\Delta C^{op}$-module (resp. $\Delta C^{op}$-module) trivial $K$ :*

$$b = \mathcal{L}_{\Sigma_\bullet}(K), \qquad \overline{b} = \mathcal{L}_{\Sigma_\bullet}(K).$$

*Démonstration :* Démontrons l'égalité dans le cas $\mathcal{F}(as)$, le cas $\Gamma(as)$ se faisant de la même manière. Nous avons vu dans la démonstration de la proposition précédente que $\Phi \in Hom_{\mathcal{F}(as)}([n],[0])$ s'écrit de manière unique sous la forme $\Phi = \Psi \circ t^k \circ \sigma$ où $\Psi$ est l'unique élément de $Hom_\Delta([n],[0])$, $0 \leqslant k \leqslant n$ et $\sigma \in \Sigma_n$. Les éléments d'une base de $b$ peuvent être vus comme étant les morphismes ensemblistes $[n] \to [0]$ ayant un ordre cyclique sur leur image réciproque. On obtient alors les relations d'équivalence suivantes sur les morphismes de $Hom_{\mathcal{F}(as)}([n],[0])$ :

$$\Psi \circ t^k \circ \sigma \sim \Psi \circ t^l \circ \sigma$$

ce qui prouve que $b = \mathcal{L}_{\Sigma_\bullet}(K)$. $\square$

**3.5 Théorème (Pirashvili-Richter [6]).** *Soient $M$ un module simplicial, $N$ un module cyclique et $K$ les modules triviaux. On a alors les isomorphismes naturels :*

$$Tor_*^{\Gamma(as)}(\overline{b}, M) = Tor_*^{\Gamma(as)}(\mathcal{L}_{\Sigma_\bullet}(K), M) \cong Tor_*^{\Delta^{op}}(K, \mathcal{O}(M))$$
$$Tor_*^{\mathcal{F}(as)}(b, N) = Tor_*^{\mathcal{F}(as)}(\mathcal{L}_{\Sigma_\bullet}(K), N) \cong Tor_*^{\Delta C^{op}}(K, \mathcal{O}(N))$$

*où $b$ et $\overline{b}$ sont les foncteurs définis dans [6]. Ces foncteurs sont en fait les images des foncteurs $K$ par les foncteurs pseudo-adjoints $\mathcal{L}_{\Sigma_\bullet}$.*

*Démonstration :* D'après la proposition 3.2, on a $\mathcal{F}(as) = \Delta C^{op} \bowtie \Sigma_\bullet$ et $\Gamma(as) = \Delta^{op} \bowtie \Sigma_\bullet$. Les isomorphismes découlent alors du corollaire 2.9 et de la proposition 3.4. $\square$

**3.6 Remarque.** Il est montré dans [1] qu'en prenant le foncteur associé à une algèbre associative $A$ et à un bimodule $M$ :

$$\begin{aligned} \mathcal{L}(A, M) : \Delta^{op} &\longrightarrow K\text{-}Mod \\ [n] &\longmapsto M \otimes A^{\otimes n}, \end{aligned}$$

on trouve l'homologie de Hochschild : $H_*(A, M) = Tor_*^{\Delta^{op}}(K, \mathcal{L}(A, M))$. De plus, on peut munir $\mathcal{L}(A, A)$ d'une structure de $\Delta C^{op}$-module, et dans ce cas on retrouve l'homologie cyclique : $HC_*(A) = Tor_*^{\Delta C^{op}}(K, \mathcal{L}(A, A))$ (voir [1]). Il est montré dans [3] que le foncteur $\mathcal{L}(A, M)$ se factorise par $\Gamma(as)$ (resp. le foncteur $\mathcal{L}(A, A)$ se factorise par $\mathcal{F}(as)$). On obtient donc les interprétations suivantes de l'homologie de Hochschild et de l'homologie cyclique :

$$\begin{cases} H_*(A, M) &\cong Tor_*^{\Gamma(as)}(\overline{b}, \mathcal{L}(A, M)), \\ HC_*(A) &\cong Tor_*^{\mathcal{F}(as)}(b, \mathcal{L}(A, A)). \end{cases}$$

# Références